\documentclass[11pt,twoside]{article}

\RequirePackage[colorlinks,citecolor=blue,urlcolor=blue]{hyperref}
\usepackage[utf8]{inputenc}
\usepackage[T1]{fontenc}
\usepackage[english]{babel}
\usepackage{charter}
\usepackage{fullpage}
\usepackage{enumerate}
\usepackage{xcolor}
\usepackage{authblk}
\usepackage{mathtools}
\usepackage{amsmath}
\usepackage{amsthm}	
\usepackage{amsfonts}
\usepackage{mathrsfs}	
\usepackage{amssymb}
\usepackage{dsfont}
\usepackage{bbm}

\theoremstyle{plain}
\newtheorem{theorem}{Theorem}[section]

\newtheorem{example}[theorem]{Example}
\newtheorem{assumption}[theorem]{Assumption}

\theoremstyle{definition}

\numberwithin{equation}{section}
\numberwithin{figure}{section}

\usepackage[nodayofweek]{datetime}

\newcommand{\bE}{\mathbb{E}}

\newcommand{\bR}{\mathbb{R}}

\newcommand{\cA}{\mathcal{A}}

\def \eps {\epsilon}

   \def\s{\sigma}

\newcommand{\bmu}{\bar \mu}

\newcommand{\cAmca}{\cA^{\textit{MCA}}}

\newcommand{\cAaea}{\cA^{\textit{AEA}}}
\newcommand{\cAcaea}{\cA^{\textit{C-AEA}}}
\newcommand{\cAesaea}{\cA^{\textit{ES-AEA}}}
\newcommand{\cAcsaea}{\cA^{\textit{CS-AEA}}}

\renewcommand{\phi}{\varphi}

\def\<{\langle}
\def\>{\rangle}

\newcommand{\law}[1]{ \mathscr L ({#1})  }
\usepackage[colorinlistoftodos,bordercolor=orange,backgroundcolor=orange!10,linecolor=orange,textsize=scriptsize]{todonotes}
\usepackage{todonotes}

\title{New particle representations for ergodic McKean-Vlasov SDEs}
\author[1]{Houssam AlRachid }
\author[2]{Mireille Bossy}
\author[3]{Cristiano Ricci}
\author[4,5,6]{Lukasz Szpruch}
\affil[1]{IDP Laboratory, Orleans University, Orleans, France}
\affil[2]{Universit{\'e} C{\^o}te d’Azur, Inria, France}
\affil[3]{University of Florence, Italy}
\affil[4]{School of Mathematics,University of Edinburgh}
\affil[5]{The Alan Turing Institute, London}
\affil[6]{Maxwell Institute for Mathematical Sciences, Edinburgh}

\date{}

\graphicspath{{pix/}}
\begin{document}	%BEGINNING
\selectlanguage{english}
\maketitle
%%%%%%%%%%%%%%%%%%%%%%%%%%%%%%%%%%%%%%%%%%%%%%%%%%%%%%%%%%%%%%%%%%%%%%%%%%%%%%%%%%%%%%%%%%%%%%%%%%%%%%%%%%%%%%%%%%%

%%%%%%%%%%%%%%%%%%%%%%%%%%%%%%%%%
\begin{abstract}
The aim of this paper is to introduce several new particle representations for \textit{ergodic} McKean-Vlasov SDEs. We construct new algorithms by leveraging recent progress in weak convergence analysis of interacting particle system. We present detailed analysis of errors and associated costs of various estimators, highlighting key differences between long-time simulations of linear (classical SDEs) versus non-linear (Mckean-Vlasov SDEs) process.   
%
%The first algorithm is an extension of the classical Euler discretization over space, where we integrate the notion of multi-clouds, sets of completely independent particles, and ergodic averages. The second one uses an approach similar to the well known Picard iteration, trying to exploit the stationary behavior and adding a multi-level discretisation to reduce the variance of the estimation. We compare numerically those two methods with respect to classical particle schemes, specifically focusing on the variance of the estimation, confronted with the error and the computational cost.
\end{abstract}

%\keywords{McKean-Vlasov equation, system of interacting particles, multi-clouds, Picard iteration, multi-level discretisation, variance reduction.}
%%%%%%%%%%%%%%%%%%%%%%%%%%%%%%%%%
{\bf 2010 AMS subject classifications:} 
Primary: 
65C30% Stochastic differential and integral equations
%,% 60H35%Computational methods for stochastic equations
; secondary: 
60H30%Applications of stochastic analysis (to PDE, etc.)
.\\
%%%%%%%%%%%%%%%%%%%%%%%%%%%%%%%%%%%%%%%%%%%%%%%%%%%%%%%%%%
\vspace{-5mm}
%\tableofcontents
%%%%%%%%%

\section{Introduction}\label{intro}

Diffusion processes are at the core of many algorithms used in statistics to sample from, typically high-dimensional, distributions. These algorithms are often based on some variant of Langevin stochastic dynamics \cite{lel10}. Given a probability measure $\pi$ (possibly known up to a normalising constant), the key idea is to construct a diffusion process which admits $\pi$ as its invariant measure. Then one can run long-time simulations of that diffusion to obtain samples from $\pi$. This ideas has been extensively studied in the context of classical SDEs \cite{talay1990second,talay2002stochastic,lamberton2002recursive,lemaire2010unconstrained,mattingly2010convergence,pages2012ergodic,dalalyan2017theoretical,durmus2017nonasymptotic}. 

 Recently new promising classes of algorithms based on the theory of gradient flows takes the form of McKean-Vlasov ODEs or SDEs \cite{csimcsekli2018sliced,bernton2018langevin,liu2017stein}. To turn them into practical algorithms one needs to approximate them with systems of interacting diffusions also called stochastic interacting particle systems. The key challenge is that, typically, with the increase of the dimension of the problem one needs to consider large number $N$ of particles. Because, for most models, the cost of particle samples growths as $N^2$ (as each particle interacts with others), the computational cost for simulating the particle systems is prohibitive. Another complication is that when using a single ensemble of particles the statistical error due to the approximation of the measure creates biased dynamics. Put differently bias is a non-linear function of the statistical error.  In addition particles are not independent. All of that renders classical variance reduction techniques not directly applicable and consequently simulations of particle systems challenging. The high computational cost is even more pronounced when the aim is to simulate particle systems over a long-time horizon. This should not come a surprise as particle systems give rise to probabilistic numerical methods for highly non-linear PDEs e.g Burgers or Navier-Stokes PDEs.   

In this work we leverage recent progress in weak convergence analysis of interacting diffusions \cite{weak-expansion,chassagneuxprobabilistic,kolokoltsov2010nonlinear,carmona2018probabilistic}. With this new insight we propose several new algorithms and analyse their errors and costs. The emphasis of the work is on algorithmic side and we gloss over some theoretical bounds that will require further research in future. As such we see this work as beacon that helps to identify the most promising research directions in the area of simulations of the ergodic  particle systems. 

\medskip

Let $(\Omega,\mathcal{F}, (\mathcal{F}_t)_{t\geq 0}, \mathbb{P})$ be a filtered probability space endowed with an  $\mathbb{R}^{k}$-valued Wiener process $w=(w_t)_{t\geq 0}$. Let $b:\bR^d \times \mathcal P (\bR^d) \rightarrow \bR^d$ and $\sigma:\bR^d \times \mathcal P (\bR^d) \rightarrow \bR^{d\times k}$. We consider, for $t\geq 0$, the McKean--Vlasov SDEs (McKV-SDE)
\begin{align}\label{eq mkvsde}	
\left\{
\begin{aligned}
&x_t = x_0 + \displaystyle \int_0^t b(x_s,\mu_s)\,ds 
+ \int_0^t \sigma(x_s,\mu_s)\,dw_s \\
& \mu_t \text{ is the law of }  x_t,
\end{aligned}\right. 
\end{align}
where $x_0$ is distributed according to a given $\mathbb{R}^{d}$-measure $\mu_0$. The  nonlinearity in the  McKean-Vlasov SDEs~\eqref{eq mkvsde} appears through the dependency of its coefficients on the law of the process. Existence of the unique solution to \eqref{eq mkvsde} has been established under various conditions on $(b,\sigma)$. See \cite{sznitman1991topics,Mele96} classical results on that topic that mainly cover the case of finite time interval. For the infinite time horizon we refer to \cite{veretennikov2006ergodic,hammersley2018mckean}.
 
Furthermore, \cite{ganz08} gives conditions for the existence  and uniqueness of the invariant measure $\pi$ for the equation \eqref{eq mkvsde}. We refer to \cite{eberle2016quantitative} for more complete theory.  In particular  \cite{eberle2016quantitative} gives fairly general conditions that guarantee that the convergence to the invariant measure in the $L^2$-Wasserstein distance  is exponentially fast for some $\lambda > 0$,  i.e 
\begin{equation} \label{eq w2 bound}
W_2(\mu_t,\pi) \leq \exp{(- \lambda t)}W_2(\mu_0,\pi).
\end{equation}
One can also control the bias of ergodic averages 
\begin{equation} \label{eq:ergo_bias}
\left| \bE\left[ \frac{1}{t} \int_0^t f(x_s) ds - \int_{\bR^d} f(x) \pi (dx) \right] \right| \lesssim  t^{-1 }.
\end{equation}
Consider the following system of $N$ particles $(x^{1,N},x^{2,N},...,x^{N,N}) $ defined as
 \begin{align}\label{eq:particlesystem}	
\left\{
\begin{aligned}
&x^{i,N}_t = x^{i,N}_0 +\displaystyle \int_0^t b(x^{i,N}_s,\mu^{N}_s )\,ds ,
+ \int_0^t \sigma(x^{i,N}_s,\mu^{N}_s)\,dw^i_s\\
&\mu^{N}_t = \frac{1}{N} \displaystyle \sum_{i=1}^N \delta_{x^{i,N}_t},
\end{aligned}\right.
\end{align}
where $(w^i)_i$ are independent $k$-dimensional Brownian motions, $(x_{0}^{i})_i$ are initial i.i.d. variables independent of $(w_t^i)_i$. The  measure valued random variable $\mu^{N}_t$ is an empirical measure of the system at time $t$. For the  purpose of computer simulations one needs to introduce time discretisation to simulate \eqref{eq:particlesystem}. We will do that in the forthcoming section. Under classical Lipschitz continuity conditions the law, seen as element of $\mathcal P([0,T],\bR^d)$, of every fixed subsystem of $k$ particles from $(x^{i,N})$ converges, when $N$ tends to infinity, to the law $\mu^{\otimes k}$. This property is called propagation of chaos phenomena. Under strong convexity of the drift, a time-uniform version of propagation of chaos has been established in \cite{ganz08}. In a  recent work, \cite{durmus2018elementary} it has been shown that in general only convexity at infinity is needed.  

Let $ f:\bR^d \rightarrow \bR$. The objective of this work is to derive, analyse and numerically investigate, several novel particle representations that will allow to approximate: 
\begin{equation} \label{eq eint}
\int_{\bR^d} f(x) \pi(dx). 
\end{equation}
To motivate our work, let's temporarily assume that  $(b,\sigma)$ do not depend on measure, i.e we are dealing with a classical SDEs. Then a typical strategy in obtaining an approximation to \eqref{eq eint} would be to set a finite time $t$ and take $N$ i.i.d. trajectories to compute $\frac{1}{t} \int_0^t \frac{1}{N} \sum_{i=1}^{N} f(x^i_s)ds$ ($N=1$ corresponds to ergodic estimator).  The error of this estimator can be decomposed as follows
\begin{align*}
&\bE\left[ \left(\int_{\bR^d} f(x) \pi (dx) - \frac{1}{t} \int_0^t \frac{1}{N}\sum_{i=1}^N\ f(x^i_s)ds \right)^2\right]^{1/2} \\
&\quad \leq  \bE\left[ \left( \int_{\bR^d} f(x) (\pi (dx) - \mu_{t}(dx))  \right)^2 \right]^{1/2}
+ \bE\left[ \left(  \int_{\bR^d} f(x)\mu_{t}(dx) - \frac{1}{t} \int_0^t \frac{1}{N}\sum_{i=1}^N\ f(x^i_s) ds \right)^2\right]^{1/2}.
\end{align*}
 The first term in the right hand side is the (weak) error of approximating the invariant measure which decays to zero 
as $\exp{(- \lambda t)}$ due to \eqref{eq w2 bound}. The second term is CLT type result and can be shown to decay to zero as $(N \cdot t)^{-1/2}$, see e.g \cite{Kipnis1986, Cattiaux2012}). We see that both $t$ and $N$ have the same impact on the variance. In the case of SDEs the cost of simulations is linear in $t$ and $N$ and hence one may be indifferent weather to simulate one long trajectory (ergodic estimator) and many shorter ones (space average). Of course if one uses parallel computer architecture, taking more samples is much more efficient. 

The situation of McKean-Vlasov SDE \eqref{eq mkvsde}	is dramatically different. The cost of simulating interacting particles \eqref{eq:particlesystem} is $N^2$ while it still increases linearly with time. As we will show it is possible to construct estimator that has one-order of magnitude lower cost while maintains the  same accuracy. Furthermore, we will investigate  \textit{ensemble} version of interacting diffusions where we generate $M$ independent systems of particles with $N$ particles in each system (ensemble). More precisely we define
\begin{align}
\label{eq cloud particle-intro}	
\left\{
\begin{aligned}
&	x^{(i,N),(j,M)}_t = x^{(i,N),(j,M)}_0 +\displaystyle \int_0^t b\big(x^{(i,N),(j,M)}_s, \mu^{N,j}_s \big)\,ds 
+ \int_0^t \sigma\big(x^{(i,N),(j,M)}_s,\mu^{N,j}_s\big)\,dw^{i,j}_{s} ,\\
& \mu^{N,j}_t  =\displaystyle \frac{1}{N}\sum_i^{N} \delta_{x^{(i,N),(j,M)}_t},
\end{aligned}\right. 
\end{align}
where $(w^{i,j},i,j)$ are independent Brownian motions. That way particles within each ensemble $j^*$ driven by $(w^{i,j^*})_{i,j^*}$ are interacting and are not independent. The particle systems $j$ and $j'$, $j \neq j'$, driven by 
$(w^{i,j})_{i,j}$ and $(w^{i,j^{'}})_{i,j}$ respectively, are independent. This idea for finite time simulations has been proposed in \cite{haji2016multilevel}. Another approach that we investigate are self-interacting diffusion  
\begin{equation} \label{eq self-interacting-intro}	
	z_t = x_0 +\displaystyle \int_0^t \left(\frac{1}{s}\int_0^s b(z_s,z_r) dr\right)ds 
+ \int_0^t \left(\frac{1}{s}\int_0^s \sigma(z_s,z_r) dr\right)dw_s 
\end{equation}
We expect that the law of $z_t$ approximates the law of $x_t$ for large $t$ due to ergodic property (\ref{eq w2 bound},\ref{eq:ergo_bias}). This gives an alternative to the particle system. We will show that the structure of the equation seems to play a crucial role in this set up.

The rest of the paper is organised as follows. In Section~\ref{setup}, we recall some classical methods to approximate \eqref{eq eint} and give error estimations and computational costs of the associated algorithms.
In Section~\ref{sec particle} we study several variants of algorithms for ergodic interacting particle systems
In Section~\ref{key idea}, we present the ensemble algorithm with ergodic average particle system. We end this paper with a general conclusion and some perspectives in Section~\ref{conc}.

\section{Setup}\label{setup}

\subsection{Algorithms}\label{algo}

As in this work we are interested in designing implementable algorithms for  \eqref{eq mkvsde}, we need to introduce time discretisation for \eqref{eq:particlesystem}. Let us define $t^n_k :=  \tfrac{k}{n}$, $k=0,1,\ldots$ and $\kappa_n(t) = t^n_k$ 
for $t\in [t^n_k, t^n_{k+1})$. We introduce (continuous time) Euler approximations for each $i\leq N$, $(y^{i,N}_t,t\geq 0)$,  $n\in \mathbb{N}$,
\begin{align}\label{eq euler}
\left\{
\begin{aligned}
&y^{i,N}_{t} = y^{i,N}_0  + \displaystyle \int_0^t b(y^{i,N}_{\kappa_n(s)},\bmu_{\kappa_n(s)}^N )\,ds 
 + \int_0^t \sigma(y^{i,N}_{\kappa_n(s)},\bmu_{\kappa_n(s)}^N)\,dw^i_{s},  \\
& \bmu_{\kappa_n(t)}^N = \displaystyle \frac{1}{N}\sum_{i=1}^N \delta_{y^{i,N}_{\kappa_n(t)}}.
\end{aligned}\right.
\end{align}

%In what follows will simply write $y^{i,N}_{t}= y^{i}_{t}$.
One may approximate \eqref{eq eint}  by ergodic average estimator with fixed $t$,
\begin{gather} \label{eq ea} \tag{EA}
\frac{1}{ t}\int_0^{ t} f(y^{1,N}_{\kappa_n(s)} ) ds\,. 
\end{gather}
For error analysis one needs to choose $N$ (number of particles) and $n$ (number of timesteps) to control the bias of the approximation of \eqref{eq mkvsde}. While \eqref{eq ea} estimator is a reasonable choice for computing approximation of \eqref{eq eint} for the invariant measures induced by classical SDEs as we already argued, the case of McKean-Vlasov SDEs approximated with a particle system, \eqref{eq ea} estimator does not seem to the best choice. This is because when using particle system \eqref{eq:particlesystem}, one computes $N$ particles and therefore calculating ergodic average along one trajectory is not efficient. See sections \ref{sec:costMCA}, \ref{sec:costAEA} for more details.  Hence, improvement can be obtained by computing, averaged ergodic average estimator  
\begin{gather} \label{eq aea} \tag{AEA}
	\frac{1}{N} \sum_{i=1}^N \left( \frac{1}{t'}\int_0^{t'} f(y^{i,N}_{\kappa_n(s)}) ds\right)\,.
\end{gather} 
Of course one expects that $t'$ in \eqref{eq aea} to be smaller than $t$ in \eqref{eq ea} for fixed accuracy. Alternative strategy for approximating \eqref{eq eint} is to resort to the standard Monte Carlo estimator where the average is taken only "over the space". More precisely we compute Monte Carlo average
\begin{gather} \label{eq mc} \tag{MCA}
\frac{1}{N} \sum_{i=1}^N f(y^{i,N}_t).
\end{gather}
Of course, the above estimator is less efficient than \eqref{eq aea}, as we explain in the coming section. The only reason we study it here is to warn practitioners that if not careful with setting up particle estimators the cost might be huge. 

\subsubsection*{Computational Cost.} 
By $\cA(\eta)$ we denote an algorithm that outputs the approximation for the quantity \eqref{eq eint}, where $\eta$ denotes the set of all the parameters we need to choose to implement it. To be able to compare the algorithms we need to fix a measure of error. For simplicity, we restore to the mean-square-error:
\begin{equation}
\textit{mse}(\cA(\eta)):= \bE \left[ \left( \int_{\bR^d} f(x) \pi (dx) - \cA(\eta) \right)^2 \right]^{1/2}.	
\end{equation}
With the measure of error of a given estimator set up, the second equally important quantity is the computation cost of algorithm $\cA$, denoted by $\textit{cost}(\cA)$. With both quantities in place 
% we can ask the question what is the optimal allocation of the parameters to achieve prescribed error tolerance. 
we can wonder about the optimal choice of parameters achieving a prescribed tolerance.
More precisely, for fixed  error tolerance $\eps>0$, we need to solve the following optimisation problem: 
\[
\left\{
\begin{aligned}
	& \textit{mse}(\cA(\eta)) < \epsilon\, ,\\
	& \min_{\eta} \textit{cost}(\cA(\eta))\,.
\end{aligned}
\right.
\]

\subsection{Assumptions}\label{sec assum}

In this section we list all the assumptions needed for our considerations. The only assumption that has not been yet established in the literature is uniform in time particle error \eqref{as weak particle error} estimation below. To the best of authors knowledge only finite time weak particle error has been studied. It is clear how to extend the weak convergence to be uniform in time but this would require lengthy introduction of heavy machinery of PDEs on measure spaces. This falls outside this paper. All other assumptions are established in literature under various level of generality and we point out reader to the corresponding papers. 

We label by $x^{i}$ a McKean-Vlasov SDE driven by $i$th Brownian motion, that is 
\begin{equation}
\label{eq:particlesystem2}
x^{i}_{t} = x^{i}_{0} + \displaystyle \int_0^t b(x^{i}_{s},\mu_s)ds + \int_0^t \s(x^{i}_{s},\mu_s)dw_s^i\,.
\end{equation}

\begin{assumption}
Convergence rate to ergodic measure: there exists $\lambda>0$ such that
\begin{gather} \label{as ergodic} \tag{HE}
 	W_2(\mu_t, \pi) \lesssim \exp{(- \lambda t)}W_2(\mu_0,\pi).
\end{gather}
\end{assumption}
As we already mentioned this has been proved in \cite{ganz08}  and \cite{eberle2016quantitative} under fairly general conditions.  
\begin{assumption}
Convergence rate of ergodic average
\begin{gather} \label{as ea} \tag{HEA}
\left( \bE \left[\left( \int_{\bR^d} f(x)\mu_{t}(dx) - \frac{1}{t} \int_0^t f(x_s)ds \right)^2\right] \right)^{1/2} \lesssim \frac{( \sup_{s\in[0,t]}\mathbb V ar[x_s])^{1/2}}{\sqrt{t}}.
\end{gather}
\end{assumption}
This is classical CLT result. See \cite{Cattiaux2012}. 
\begin{assumption}
Uniform in time weak convergence of the particle system: for sufficiently smooth~$f$, 
\begin{gather} \label{as weak particle error} \tag{HW}
 \sup_{t\geq 0}| \bE f(x^1_t) - \bE f(x^{1,N}_t) | \lesssim \frac{1}{N}.
\end{gather}
\end{assumption}
This type of bound is new in the literature. We refer a reader to \cite{weak-expansion,chassagneuxprobabilistic,kolokoltsov2010nonlinear,carmona2018probabilistic} for more details. 
\begin{assumption}
Uniform in time strong propagation of chaos
\begin{gather} \label{as strong particle error} \tag{HS}
\sup_{t\geq 0}(\bE| x^1_t -  x^{1,N}_t |^2)^{1/2}  \lesssim \frac{1}{\sqrt{N}}.
\end{gather}
\end{assumption}
See \cite{durmus2018elementary} for details.
\begin{assumption}
Uniform in time weak discretisation error: for sufficiently smooth $f$, 
\begin{gather} \label{as weak discretisation error} \tag{HDW}
 \sup_{t\geq 0}| \bE f(x^1_t) - \bE f(y^{1,N}_t) | \lesssim \frac{1}{n}.
\end{gather} 
Uniform in time strong discretisation error
 \begin{gather} \label{as strong discretisation error} \tag{HDS}
\sup_{t\geq 0}(\bE| x^1_t -  y^{1,N}_t |^2)^{1/2}  \lesssim \frac{1}{\sqrt{n}}.
\end{gather}
\end{assumption}
One can refer to \cite{ganz08}, where such results are proved.

%\ls{Remark on the assumptions: \ref{as ergodic} holds for strongly convex potential (e.g Bossy, Talay, Ganz); \ref{as weak particle error} has been recenlty established as by product of Ito formula for measures (apparently there is a comment on that in the book by Carmona, Delarue, but non of these results is uniform in time; \ref{as strong particle error} is classical when both coefficients are Lipschitz, Sznitman for finite time interval, moreover in papers on ``uniform propagation of chaos'' under strong convexity of the the drift it is rather simple excerice; weak and strong errors in this setup are classical (again for Lipschitz) and one can refer to PhD thesis by Bossy, uniform in time   Bossy, Talay, Ganz ?? }
\section{Algorithms for Ergodic Interacting Particle systems}\label{sec particle}

\subsection{Monte Carlo Average}\label{sec:costMCA}
%\ls{It is important to decompose the error in the correct way in order to take full advantage from the assumptions stated in the previous section.}
For the Monte Carlo Average estimator we introduce  the following notation
$\cAmca(t,n,N)$. The aim is to find the optimal allocation of the parameters $ (t,n,N)$ for fixed mean-square-error. We have
\begin{align*}
	\textit{mse}(\cAmca(t,n,N)) 
	=  & \bE \left[ \left( \int_{\bR^d} f(x)\pi(dx) - \frac{1}{N} \sum_{i=1}^N f(y^{i,N}_t)  \right)^2 \right]^{1/2}	\\
& \quad \lesssim     \left| \int_{\bR^d} f(x) (\pi (dx) - \mu_{t}(dx))  \right|  +  \left| \bE[ f(x_{t}) ]  - \bE[f(x^{1,N}_{t})]  \right|     \\
& \quad \quad +  \left| \bE[  f(x^{1,N}_{t})]    - \bE[ f(y^{1,N}_t)]  \right|   +  \bE \left[ \left( \bE[ f(y^{1,N}_t)]   - \frac{1}{N} \sum_{i=1}^N f(y^{i,N}_t)  \right)^2 \right]^{1/2} .	
\end{align*}
The four error terms are in order: bias (due to finite time simulation); weak particle approximation error; weak time discretisation error; variance/propagation of chaos. 
The first three error terms can be estimated directly from the Assumptions in Section \ref{sec assum}. The last variance error term requires extra comment, as the fact that particles are not i.i.d does not allow to use classical central limit theorem (CLT). Indeed
\begin{equation*}
\begin{split}
	 \bE [ (  \bE[ f(y^{1,N}_t)]   - \frac{1}{N} \sum_{i=1}^N f(y^{i,N}_t))^2 ]^{1/2} 
	 = &\, \bE [ ( \bE[\frac{1}{N} \sum_{i=1}^N f(y^{i,N}_t)]   - \frac{1}{N} \sum_{i=1}^N f(y^{i,N}_t))^2]^{1/2} \\
	 \leq & \, \bE[  ( \frac{1}{N}\sum_{i=1}^N f(y^{i,N}_t) )^2]^{{1/2}}.
\end{split}	
\end{equation*}
Next, define $(\tilde{y}^{i})_i$ as the solution of the continuous Euler scheme:
\begin{equation*}
\begin{cases}
			\tilde{y}^{i}_{t} = \tilde{y}^{i}_0  + \displaystyle \int_0^t b(\tilde{y}^{i}_{\kappa_n(s)},\law {\tilde{y}_{\kappa_n(s)}}  )\,ds 
 + \int_0^t \sigma(\tilde{y}^{i}_{\kappa_n(s)},\law{\tilde{y}_{\kappa_n(s)}} )\,dw^i_{s}\,,  \\
 \law{\tilde{y}_{t}}   =\textit{Law} (\tilde{y}_{t})\,.
\end{cases}
\end{equation*}
It is an easy exercise to show that strong propagation of chaos (HS) can be established on the level of Euler discretisation. This together with Cauchy-Schwarz inequality gives 
\begin{equation}\label{eq variance}
\begin{split}
&\bE \left[  \left(  \frac{1}{N}\sum_{i=1}^N f(y^{i,N}_t) -  f(\tilde{y}^{i}_t) \right)^2\right] \\
&= \bE \left[   \frac{1}{N^2}\sum_{i=1}^N (f(y^{i,N}_t) -  f(\tilde{y}^{i}_t))^2   + \frac{1}{N^2}\sum_{i<j}^N (f(y^{i,N}_t) -  f(\tilde{y}^{i}_t) )f(y^{j,N}_t) -  f(\tilde{y}^{j}_t) )   \right] \\
& \leq    \frac{1}{N^2}\sum_{i=1}^N \bE\left[(f(y^{i,N}_t) -  f(\tilde{y}^{i}_t) )^2\right]  \\
&\qquad\qquad + \frac{1}{N^2}\sum_{i<j}^N (\bE[(f(y^{i,N}) -  f(\tilde{y}^{i,N}_t) )^2])^{1/2}(\bE[(f(y^{j}_t ) -  f(\tilde{y}^{j}_t) )^2])^{1/2} \lesssim \frac{1}{N}\,.
\end{split} 
\end{equation}
This, and the fact that $(\tilde{y}^{i}_{t})_i$ are i.i.d. allows to conclude that
\begin{equation*}
\begin{split}
\bE\left[  \left( \frac{1}{N}\sum_{i=1}^N f(y^{i,N}_t ) \right)^2\right]
 = & 
\, \bE \left[  \left( \frac{1}{N}\sum_{i=1}^N f(\tilde{y}^{i}_t ) +  ( f(y^{i,N}_t ) -  f(\tilde{y}^{i}_t ) \right)^2\right] \\
\lesssim & \,
\bE \left[  \left( \frac{1}{N}\sum_{i=1}^N f(\tilde{y}^{i}_t ) \right)^2\right]  
  + \bE \left[  \left(  \frac{1}{N}\sum_{i=1}^N f(y^{i,N}_t ) -  f(\tilde{y}^{i}_t ) \right)^2\right]
  \lesssim  \frac{1}{N}.
\end{split}	
\end{equation*}

%Next, define $\t ilde{y}$ as the solution of the continuous Euler scheme:
%\begin{equation*}
%\begin{cases}
%			\tilde{y}^{i}_{t} &= \tilde{y}^{i}_0  + \displaystyle \int_0^t b[\tilde{y}^{i}_{\kappa_n(s)},\law {\tilde{y}_{\kappa_n(s)}} ^N ]\,ds 
% + \int_0^t \sigma[\tilde{y}^{i}_{\kappa_n(s)},\law{\tilde{y}_{\kappa_n(s)}} ^N]\,dw_{s}\,,  \\
% \law{\tilde{y}_{t}} ^N  &=Law  (\tilde{y}_{t})\,.
%\end{cases}
%\end{equation*}
%Using the fact that $(\tilde{y}^{i}_{t})_i$ are i.i.d and applying the  strong propagation of chaos property for Euler scheme ($f$ is Lipschitz), then:
%\begin{equation}
%\begin{split}
%\bE\left[  \left( \frac{1}{N}\sum_{i=1}^N f(y^{i,N}_t ) \right)^2\right]
% = & 
%\, \bE \left[  \left( \frac{1}{N}\sum_{i=1}^N f(\tilde{y}^{i}_t ) +  ( f(y^{i,N}_t ) -  f(\tilde{y}^{i}_t ) \right)^2\right] \\
%\lesssim & \,
%\bE \left[  \left( \frac{1}{N}\sum_{i=1}^N f(\tilde{y}^{i}_t ) \right)^2\right]  
%  + \bE \left[  \left(  \frac{1}{N}\sum_{i=1}^N f(y^{i,N}_t ) -  f(\tilde{y}^{i}_t ) \right)^2\right] \\
%  \lesssim & \frac{1}{N}.
%\end{split}	
%\end{equation}
From here and Assumptions in Section \ref{sec assum}  we have
$$
\textit{mse}(\cAmca(t,n,N) )\lesssim e^{{-}\lambda t} + \frac{1}{N} + \frac{1}{n} +\frac{1}{\sqrt{N}}.
$$
Notice that because of the  term $1/\sqrt{N}$, it is not clear how we can take advantage of the assumption \eqref{as weak particle error}. Fix $\epsilon >0$ and set 
$
\textit{mse}(\cAmca(t,n,N)) \lesssim \epsilon.
$
This leads to the following choice of the parameters 
$
t \approx \lambda^{-1} \log(\epsilon^{-1}), N \approx \epsilon^{-2},  n \approx {\epsilon}^{-1}. 
$
As the cost of simulating particle system at every step of the Euler scheme is $N^2$ we have 
$$
\textit{cost}(\cAmca(t,n,N)) = t n N^2 \approx \log(\epsilon^{-1}) \epsilon^{-5}. 
$$ 
This should be compared to $tnN = \log(\epsilon^{-1}) \epsilon^{-3}$ for the simulation of standard SDEs. 
%while in the linear 
%\[
%\textit{cost}(\cAmca(t,n,N)) = t n N \approx \log(\epsilon^{-1}) \epsilon^{-3}.
%\]  

\subsection{Averaged Ergodic Average}\label{sec:costAEA}
As before, we denote by $\cAaea(t,n,N)$ the averaged ergodic average estimator in \eqref{eq aea}. We have
\begin{align*}
	\textit{mse}(\cAaea(t,n,N)) = 
	 & \bE \left[ \left( \int_{\bR^d} f(x) \pi (dx) 
	 - \frac{1}{N} \sum_{i=1}^N \left( \frac{1}{t}\int_0^{t} f(y^{i,N}_{\kappa_n(s)}) ds\right) \right)^2 \right]^{1/2}	\\
 &   \lesssim \ \left| \int_{\bR^d} f(x) (\pi(dx) - \mu_{t}(dx))  \right|  	
 +  \left| \bE[ f(x_{t}) ]  - \bE[f(x^{1,N}_{t})]  \right|      \\
& \quad+  \left| \bE[  f(x^{1,N}_{t})]    - \bE[ f(y^{1,N}_t)]  \right|   
+  \bE \left[ \left( \bE[ f(y^{1,N}_t)]   - \frac{1}{N} \sum_{i=1}^N \left( \frac{1}{t}\int_0^{t} f(y^{i,N}_{\kappa_n(s)}) ds\right) \right)^2 \right]^{1/2}.	
\end{align*}
To estimate the variance term we note that
\[
 \bE[ f(y^{1,N}_t)]   - \frac{1}{N} \sum_{i=1}^N \left( \frac{1}{t}\int_0^{t} f(y^{i,N}_{\kappa_n(s)}) ds\right) 
=  \bE[ \frac{1}{N}\sum_{i=1}^N(y^{i,N}_t)]   - \left( \frac{1}{t}\int_0^{t} \frac{1}{N} \sum_{i=1}^N  f(y^{i,N}_{\kappa_n(s)}) ds\right)  .	
\]
Hence to estimate the variance we use \eqref{as ea} combined with the computation \eqref{eq variance}. 
 Therefore, by the Assumptions in Section  \ref{sec assum}, we have
\[
\textit{mse}(\cAaea(t,n,N))\lesssim e^{{-}\lambda t} + \frac{1}{N}  + \frac{1}{n} + \frac{1}{\sqrt{tN}}.
\]
We notice that comparing to the \eqref{eq mc} case, the last term is multiplied by 
$1/\sqrt{t}$. The (asymptotic) cost of the algorithm is the same as before.
% and while for the linear interacting case the above algorithm is more efficient only by term $\log(\epsilon^{-1})$, for the nonlinear interacting case the difference is dramatic as the cost in time is linear while in number of particles quadratic. 
Again we fix $\epsilon$. The following choice of the parameters ensures that $\textit{mse}(\cAaea(t,n,N)) \lesssim \epsilon^2$ is 
$t \approx  \epsilon^{-1}$, $N \approx \epsilon^{-1}$, $n \approx \epsilon^{-1}$.
The cost consists of two parts: the cost of simulating particle system and the cost of computing averaged ergodic estimator. We have
\[
\textit{cost}(\cAaea(t,n,N))= t n N^2 +tN\approx \epsilon^{-4}\,.
\]
Which is an order of magnitude lower than for Monte Carlo average! 

Notice that similar computation for ergodic average estimator gives $\textit{mse}(\cA^{\textit{EA}}(t,n,N))\lesssim e^{{-}\lambda t} + \frac{1}{N}  + \frac{1}{n} + \frac{1}{\sqrt{t}}$,  leading to the same cost than Monte Carlo average. 
%{while in the linear case
%\[
%\textit{cost}(\cAaea(t,n,N))= t n N +tN\approx \epsilon^{-3}. 
%\]  }
%The improved complexity is possible due to weak error in particle approximation stated in the Assumption \eqref{as weak particle error}.

\subsection{ensemble AEA}

For the ensemble version of the algorithm we generate $M$ independent systems of particles with $N$ particles in each system. More precisely we define, 
\begin{align}
\label{eq cloud particle}	
\left\{
\begin{aligned}
&	x^{(i,N),(j,M)}_t = x^{(i,N),(j,M)}_0 +\displaystyle \int_0^t b(x^{(i,N),(j,M)}_s, \mu^{N,j}_s )\,ds 
+ \int_0^t \sigma(x^{(i,N),(j,M)}_s,\mu^{N,j}_s)\,dw^{i,j}_{s} ,\\
& \mu^{N,j}_t  =\displaystyle \frac{1}{N}\sum_i^{N} \delta_{x^{(i,N),(j,M)}_t},
\end{aligned}\right. 
\end{align}
where $(w^{i,j},i,j)$ are independent Brownian motions. That way particle within each cloud $j^*$ driven by $(w^{i,j^*})_{i,j^*}$ are interacting and are not independent. The particle systems $j$ and $j'$, $j \neq j'$, driven by 
$(w^{i,j})_{i,j}$ and $(w^{i,j^{'}})_{i,j}$ respectively, are independent. This idea has been proposed in \cite{haji2016multilevel} for the finite time simulations. 
\medskip

We consider ensemble version of  \ref{eq aea}. 
\begin{gather} \label{eq c-aea} \tag{C-AEA}
	\frac{1}{M} \sum_{j=1}^M\frac{1}{N}\sum_{i=1}^N \left( \frac{1}{t}\int_0^{t} f(y^{(i,N),(j,M)}_{\kappa_n(s)}) ds\right)\,.
\end{gather} 
In fact, all algorithms that we study can have their ensemble versions. By denoting $\cAcaea$  the new method 
\begin{align*}
& \textit{mse}(\cAcaea(t,n,N,M)) = 
\bE \left[ \left( \int_{\bR^d} f(x) \pi (dx) 
- \frac{1}{M} \sum_{j=1}^M\frac{1}{N} \sum_{i=1}^N \left( \frac{1}{t}\int_0^{t} f(y^{(i,N),(j,N)}_{\kappa_n(s)}) ds\right) \right)^2 \right]^{1/2}\\
& \qquad \lesssim     \left| \int_{\bR^d} f(x) (\pi (dx) - \mu_{t}(dx))  \right|    
 +   \left| \bE[ f(x_{t}) ]  - \bE[f(x^{1,N}_{t})]  \right|      \\
& \qquad +  \left| \bE[  f(x^{1,N}_{t})]    - \bE[ f(y^{1,N}_t)]  \right|   
 +  \bE \left[ \left( \bE[ f(y^{1,N}_t)]   - 
\frac{1}{M} \sum_{j=1}^M \frac{1}{N} \sum_{i=1}^N \left( \frac{1}{t}\int_0^{t} f(y^{(i,N),(j,M)}_{\kappa_n(s)}) ds\right) \right)^2 \right]^{1/2}.
\end{align*}
The Assumptions in Section \ref{sec assum} yield
$$
\textit{mse}(\cAcaea(t,n,N,M) )\lesssim e^{-\lambda t} + 1/N + 1/n+ \frac{1}{\sqrt{tNM}}.
$$
We notice that comparing to the previous case the last term is multiplied by 
$1/{\sqrt{M}}$. Crucially, the cost of the algorithm growths linearly in $M$. As the cost growths as $N^2$, we are better off taking $M \approx N$ to balance the error in the last term (instead of taking $M=1$ and $N^2$). To make it precise  we fix $\epsilon$. The following choice of the parameters ensures that $\textit{mse}(\cAcaea(t,n,N)) \lesssim~\epsilon$.
$$
\mbox{ $t \approx \lambda^{-1} \log(\epsilon^{-1})$, $N \approx \epsilon^{-1}$, $n \approx \epsilon^{-1}$, $M=(\lambda^{-1} \log(\epsilon^{-1}))^{-1} \epsilon^{-1}$.}$$
The cost of simulating particles and computing the estimator is
$$
\textit{cost}(\cAcaea(t,n,N)) = t n N^2 M +tNM \approx \epsilon^{-4}.
$$
This is the same as for averaged ergodic estimator. However the above computations do not take under consideration the fact that ensemble algorithms can take full advantage from the parallel computer architecture and therefore will be superior in practice.     

\section{Algorithms for Sef-interacting Particle systems}\label{key idea}

In this section, we present the key ideas improvement in the definition of more efficient algorithm.  From the decomposition of the mean square error, we see that different algorithms that we considered only affected the ``variance'' of the final estimator. Therefore, to improve the efficiency of the algorithm we need to either modified particle system itself or consider different simulation strategies such as Multilevel-Monte Carlo.  Here we focus on the former.

%\subsection{Multi-clouds algorithm with Ergodic Average particle system}
Significance of the ergodic theorem is that one can approximate the integral \eqref{eq eint} by simulating only one path of the process \eqref{eq mkvsde} rather then the whole particle system.  From now on, we will keep the structural assumptions on the coefficients  of \eqref{eq mkvsde}, namely 
\begin{align*}
\left\{
\begin{aligned}
&x_t = x_0 + \displaystyle \int_0^t b(x_s,\law {x_s})\,ds 
+ \int_0^t \sigma(x_s, \law {x_s})\,dw_s \\
&\law {x_t} \text{ is the law of }  x_t.
\end{aligned}\right. 
\end{align*}
One may consider 
\begin{equation} \label{eq self-interacting}	
	z_t = x_0 +\displaystyle \int_0^t \left(\frac{1}{s}\int_0^s b(z_s,z_r) dr\right)ds 
+ \int_0^t \left(\frac{1}{s}\int_0^s \sigma(z_s,z_r) dr\right)dw_s \,.
\end{equation}
We expect that the law of $z_t$ approximates the law of $x_t$ for large $t$ due to ergodic property \eqref{eq w2 bound}

Processes of the form  \eqref{eq self-interacting} are known in literature as self interacting diffusions. We refer to  \cite{Kleptsyn2012}  where the convergence to the invariant measure has been established.  Notice that  there is no need for the particle system any more as one could simply simulate one path of the process to calculate ergodic integral \eqref{eq eint}.

However, motivated by computations in the previous section where mixed ergodic/Monte Carlo average we introduce the corresponding mean self-integrated SDE
\begin{equation}
\label{eq self particle}	
z_t = x_0 + \int_0^t \left(\frac{1}{s}\int_0^s b(z_s,\law{z_r}) dr\right)ds 
+\int_0^t \left(\frac{1}{s}\int_0^s \bE \sigma(z_s,\law{z_r}) dr \right)dw_s 
\end{equation} 
and its independent copies $(z^i)$ driven by Brownian motion $(w^i)$. Note that "one-particle" approximation of \eqref{eq self particle}	 is precisely a self-interacting diffusion. 

To gain better insight, into the idea of using self-interacting diffusions to approximate McKean-Vlasov SDEs we consider a simple example first.

\begin{example}\label{example1}
Consider a simple scalar McKean-Vlasov SDE $x$, together with its mean self-integrated version $y$, and its self-interaction motion SDE $z$, 
\begin{align*}
&x_t = x_0  -\int_0^t \alpha x_s ds + \int_0^t \beta\bE[x_s] ds + w_t,\\
&z_t = x_0  - \int_0^t \alpha z_s ds + \int_0^t \frac{1}{s}\left(\int_0^s\beta \bE[z_\theta] d\theta\right) ds + w_t.
\end{align*}
We stress out that dissipativity comes from the part of the drift that do not depend on measure. 
We assume  $\alpha > \beta$. 
To estimate the convergence rate to the invariant measure we analyse the evolution of the difference of two solutions to the above $x$ SDE initiated at $L^2$ random variables $\xi_1$ and $\xi_2$. With
\begin{equation}\label{eq:aux1}
		\bE[ (x_t^{\xi_1} - x_t^{\xi_2})]= e^{-(\alpha - \beta) t}
		\bE [(\xi_1 - \xi_2)], 
\end{equation}
we have
\begin{equation*}
	\begin{split}
	e^{2 \alpha t}	\bE[ (x_t^{\xi_1} - x_t^{\xi_2})^2 ]= &
		\bE [(\xi_1 - \xi_2)^2]
		+ 2 \beta \int_0^t e^{2\alpha s}  (\bE[x_s^{\xi_1} - x_s^{\xi_2}] )^2ds \\
		=  & \bE [(\xi_1 - \xi_2)^2]
		+ 2 (\bE [\xi_1 - \xi_2])^2\beta \int_0^t e^{2\alpha s} e^{-2 (\alpha - \beta) s }  ds \\	
		=  & \bE [(\xi_1 - \xi_2)^2]
		+  (\bE [\xi_1 - \xi_2])^2 (e^{2 \beta t }-1).
	\end{split}
\end{equation*}
Due to properties of $W_2$ distance and the fact that the above calculation does not depend on a particular choice of random variables  $\xi_1$ and $\xi_2$, we have 
\[
W_2^2( \mathscr L(x_t^{\xi_1}),\mathscr L ( x_t^{\xi_2}) )		\lesssim e^{- 2(\alpha -\beta) t} W_2^2(\mathscr L(\xi_1), \mathscr L (\xi_2))\,.
\]
Furthermore, in this simple example we can take advantage from the explicit solutions to calculate 
\begin{align*}
	\bE[x_t] - \frac{1}{t}\int_0^t\bE[x_s]ds 
	= \bE[x_0]\left( e^{-(\alpha - \beta)t} - \frac{1}{t(\alpha - \beta)}(1-e^{-(\alpha- \beta)t}) \right)\,.
\end{align*}
Hence if $\log t < (\alpha - \beta)t$ then we have that 
\begin{align*}
	| \bE[x_t] - \frac{1}{t}\int_0^t\bE[x_s]ds | \leq c 1/t\,.
\end{align*}
\medskip
Let us consider now process $z$. By integration by part we have that 
\begin{align*}
\bE[z_t] =  \bE[z_0] -\int_0^t \alpha \bE[z_\theta]d\theta + \beta\int_0^t \log(t/\theta) \bE z_\theta d\theta
\end{align*}
and 
\begin{align*}
e^{\alpha t} \bE[z_t] = \bE[z_0] +\int_0^t e^{\alpha s}\frac{1}{s}\int_0^s \bE z_\theta d\theta\,.
\end{align*}
Then we observe that if $\bE[z_0] > 0$, then $\bE[z_t]$ stays non-negative  and do not cross 0. We observe also that $t\mapsto \bE[z_t]$ is decreasing as for all $t\geq 0$, for all  $s\in[te^{-\alpha/\beta }, t]$
\begin{align*}
\bE[z_t] - \bE[z_s] = \int_s^t \left(\beta \log(t/\theta) -\alpha\right) \bE z_\theta \,d\theta \leq 0. 
\end{align*}
In particular 
\begin{align*}
\frac{\bE[z_t] - \bE[z_s]}{t-s} = \frac{1}{t-s} \int_s^t \left(\beta \log(t/\theta) -\alpha\right) \bE z_\theta d\theta \leq  \frac{1}{t-s} \int_s^t \left(\beta \log(t/\theta) -\alpha\right) d\theta. 
\end{align*}
Taking limit $s\rightarrow t$, we obtain after integration that 
\begin{align}\label{eq:aux2}
\bE[z_t] \leq \bE[z_0] \exp(- (\alpha -\beta) t)\,.
\end{align}

The following computation on $z$  is a tentative to evidence the rate's gain in the  convergence rate to equilibrium. To this aim, we analyse the evolution of the difference of two solutions $z$ initiated at $L^2$ random variables $\xi_1$ and $\xi_2$.

Repeating the previous computation for $\bE[\xi_1-\xi_2] \geq 0$, we also obtain
\begin{equation*}
0\leq \bE[ (z_t^{\xi_1} - z_t^{\xi_2})]\leq  e^{-(\alpha - \beta) t} \bE [(\xi_1 - \xi_2)], 
\end{equation*}
and we use this  weak estimation to derive  the  $L^2$-norm bound. We have first that  
\begin{align*}
 e^{2\alpha t} \bE [ (z_t^{\xi_1} - z_t^{\xi_2})^2]  & =  
\bE [({\xi_1} -{\xi_2})^2]
		+ 2{\beta}   \int_0^t e^{2\alpha s}  \bE[z_s^{\xi_1} - z_s^{\xi_2}] \frac{1}{s}\left(\int_0^s\bE[z_\theta^{\xi_1}-z_\theta^{\xi_2}] d\theta \right)ds \\
& \quad \leq 
\bE [({\xi_1} - {\xi_2})^2]
		+ 2{\beta}   (\bE [\xi_1 - \xi_2])^2  \int_0^te^{2\alpha s}   e^{-(\alpha -\beta)  s}\frac{1}{(\alpha-\beta)s} (1 - e^{-(\alpha -\beta) s})ds\\
		& \quad \leq  
\bE [({\xi_1} - {\xi_2})^2]
		+ 2{\beta}   (\bE [\xi_1 - \xi_2])^2  \int_0^t\frac{1}{(\alpha-\beta)s} (e^{(\alpha+\beta) s } - e^{ 2\beta s})ds\,.
\end{align*}
As $\alpha >\beta$, it is not difficult to check  that  for $t$ big enough, 
\begin{align*}
0\leq &\int_0^t \frac{1}{s}(e^{ (\alpha+\beta) s} - e^{2\beta s})ds \leq  \frac{2}{t} \left(\frac{e^{(\alpha+\beta)t}}{\alpha+\beta} - \frac{e^{2\beta t}}{2\beta}\right) - \frac{2}{t}\left( \frac{1}{\alpha+\beta} - \frac{1}{2\beta}\right)\,.
\end{align*}
And we obtain the contraction inequality, 
\begin{align*}
 \bE [ (z_t^{\xi_1} - z_t^{\xi_2})^2]  & \leq 
 e^{-2\alpha t } 
		\bE [({\xi_1} - {\xi_2})^2]
		+ c \frac{e^{-(\alpha-\beta) t}}{t}  (\bE [\xi_1 - \xi_2])^2, 
\end{align*}
which means that after a time $t\geq t_0$, the convergence in $L^2$ is exponentially fast with a rate $\alpha \wedge (\alpha+ \log(t)-\beta)/2$, that accelerates and becomes with time  better than the rate for process $x$ (in $\alpha - \beta$).
\medskip

We have also to show that $x$ and $z$ have the same equilibrium measure. A sufficient condition is the $L^2$-convergence for $z$ toward $x$ in time. To not spend to much time on this example, we make the following assumption from the previous Wasserstein contraction for $z$ and from what we obtain for $x$,  that 
\begin{align*}
	| \bE[z_t] - \frac{1}{t}\int_0^t\bE[z_s]ds | \leq c 1/t\,.
\end{align*} 
We consider now
\begin{align*}
e^{2(\alpha -\beta) t}  \bE (x_t - z_t)^2 & \leq e^{2(\alpha - \beta) }\bE (x_1 - z_1)^2 - 2\beta \int_1^t  e^{2(\alpha - \beta) s}  \bE (x_s - z_s)^2 ds  \\
& \quad  + 2 \beta \int_1^t e^{2(\alpha - \beta) s} \bE (x_s - z_s) \left(
\bE[x_s]  - \frac{1}{s}\int_0^s \bE[z_\theta] 
d\theta \right)ds.
\end{align*}
But since $|\bE (x_s - z_s)| \leq |\bE (x_s)| + |\bE(z_s)| \leq c e^{-(\alpha -\beta) s}$ from \eqref{eq:aux1} and \eqref{eq:aux2}, 
\begin{align*}
& \int_1^t e^{2(\alpha  - \beta) s} \bE (x_s - z_s) \left(
\bE[x_s]  - \frac{1}{s}\int_0^s \bE[z_\theta]  
d\theta \right)ds \\
& \leq \int_1^t e^{2(\alpha -\beta) s}\bE (x_s - z_s)^2 + 
 \int_1^t e^{2(\alpha-\beta) s} |\bE (x_s - z_s)| \frac{c}{s}  ds \leq \int_1^t e^{2(\alpha -\beta) s}\bE (x_s - z_s)^2 + 
 \int_1^t e^{(\alpha-\beta) s} \frac{c'}{s}  ds 
\end{align*}
and we have obtained 
\begin{align*}
e^{2(\alpha -\beta) t}  \bE (x_t - z_t)^2  \leq e^{2(\alpha - \beta) }\bE (x_1 - z_1)^2 
+ 2\beta  \int_1^t e^{(\alpha-\beta) s} \frac{c'}{s}  ds \\
 \bE (x_t - z_t)^2  \leq e^{-2(\alpha - \beta)(t-1) }\bE (x_1 - z_1)^2 + C e^{- (\alpha -\beta) t} \log(t)\,.
\end{align*}
The convergence is then ensured. Hence, mean self-interacting version $z$ of $x$ is converging to $x$ in $L^2$, which means that $z$ convergence  to the equilibrium measure $\pi$ and hence can used as an alternative model for sampling. 
\end{example}

\subsubsection*{Self-averaged ergodic averaged algorithm}
By considering error decomposition studies in the previous section, we see that the key ingredient that we ought to understand is the following rate of convergence

Key ideas: here we want to take part of a special structure, i.e
\begin{align*}
\left\{
\begin{aligned}
&x_t = x_0 + \displaystyle \int_0^t V(x_s )  + W(x_s,\law {x_s})\,ds 
+ \int_0^t \sigma\,dw_s \\
&\law {x_t} \text{ is the law of }  x_t,
\end{aligned}\right. 
\end{align*}
such that potential $V$ is say convex and $W$ has a small Lipschitz constant we should obtain exponential "forgetting property", as we observed in Example \ref{example1}. 

For small time this will be bad approximation and because we average the error the bad approximation from the initial time will prevail. 

Further, we consider particle system of the form 
 \begin{equation}
\label{eq self particle}	
	z^{i,N}_t = x^{i,N}_0 + \int_0^t \frac{1}{N} \sum_{j=1}^N \left(\frac{1}{s}\int_0^s b(z^{i,N}_s,z^{j,N}_r) dr\right)ds 
+\int_0^t \frac{1}{N} \sum_{j=1}^N  \left(\frac{1}{s}\int_0^s \sigma(z^{i,N}_s,z^{j,N}_r) dr \right)dw^i_s 
\end{equation} 
for which the error to be investigate is
\[
\left| \bE\left[ f(z^{1}_{t}) \right]  - \bE\left[f(z^{1,N}_{t})\right]\right |, 
\] 
to have in place of \eqref{as weak particle error}. According to our computation on {example} \ref{example1}, we impose the following assumption

%\ls{Correction from here plus new complexity }
The following bound gives a leading error term. We chose $\lambda$ as a exponent for simplicity as it particular value does not affect asymptotic cost/error analysis. 
\begin{assumption}
\begin{gather} \label{as weak ergodic particle error} \tag{HEW}
 \left| \bE f(z^{1}_t) - \bE[f(z^{1,N_t}_{t})] \right| \lesssim \frac{e^{-\lambda t}}{N_t }\,.
\end{gather}
\end{assumption}

The time-discretisation of the equation \eqref{eq self particle}  reads as 	
 \begin{equation}
\label{eq self particle discrete-time}	
r^{i,{N_t}}_t = x^{i,{N_t}}_0 + \int_0^t \frac{1}{{N_t}} \sum_{j=1}^{N_t} \left(\frac{1}{s}\int_0^s b(r^{i,{N_t}}_{\kappa_n(s)},r^{j,{N_t}}_{\kappa_n(\theta)}) d\theta\right)ds 
+\int_0^t \frac{1}{{N_t}} \sum_{j=1}^{N_t} \left( \frac{1}{s}\int_0^s \sigma(r^{i,{N_t}}_{\kappa_n(s)},r^{j,{N_t}}_{\kappa_n(\theta)}) d\theta \right)dw^i_{s}.  
\end{equation} 
Let's introduce the corresponding estimator, an ensemble self-integrated version of  \ref{eq aea}. 
\begin{gather} \label{eq sc-aea} \tag{CS-AEA}
\frac{1}{M}\sum_{j=1}^{M}\frac{1}{{N_t}} \sum_{i=1}^{N_t} \left( \frac{1}{t}\int_0^{t} f(r^{(i,{N_t}),(j,M)}_{\kappa_n(s)}) ds\right)\,.
\end{gather} 
In what follows, we test \eqref{as weak ergodic particle error} in the cost analysis.  Notice that we cannot test it directly, as we need to do Monte-Carlo approximation for the expectation. We first observe that by exchangeability of the law of the particle systems, 
\begin{align*}
\bE \left[f(z^1_t)\right] - \bE\left[f(z^{1,N}_{t})\right] = 
\bE \left[ \frac{1}{N} \sum_{i=1}^{N} f(z^i_t)\right] - \bE\left[ \frac{1}{N} \sum_{i=1}^{N}  f(z^{i,N}_{t}) \right]\,.
\end{align*}
%We use this representation to decrease the variance. 
Next we consider  $M$ independent ensembles to approximate expectations that is
\begin{align*}
	& \left|  \frac{1}{M} \sum_{j=1}^{M }\frac{1}{N} \sum_{i=1}^{N} f(z^{ij}_t) - \frac{1}{M} \sum_{j=1}^{M } \frac{1}{N} \sum_{i=1}^{N}  f(z^{(i,N),(j,M)}_{t})] \right| \\& \quad \leq \left|\frac{1}{M} \sum_{j=1}^{M }\frac{1}{N} \sum_{i=1}^{N} f(z^{ij}_t) - \bE[ f(z^{1}_{t}) ]\right| \\
	 &\quad \quad +\left|  \bE[ f(z^{1}_{t}) ] - \bE[f(z^{1,N}_{t})] \right| + \left|  \bE[f(z^{1,N}_{t})]- \frac{1}{M} \sum_{j=1}^{M } \frac{1}{N} \sum_{i=1}^{N}  f(z^{(i,N),(j,M)}_{t}) \right|  \\
	 & \quad \lesssim (\sqrt{M N})^{-1} + (Nt)^{-1} + (\sqrt{M N})^{-1}
	\end{align*}
where the first and last error are standard MC estimates while the middle one is given by \eqref{as weak ergodic particle error}.

\subsection{Cost analysis}
We analyse the cost of the self-averaged ergodic averaged estimator 

\begin{align*}
	\textit{mse}(\cAesaea(t,n,N)) = 
	 & \bE \left[ \left( \int_{\bR^d} f(x) \pi (dx) 
 - \frac{1}{{N_t}} \sum_{i=1}^{N_t} \left( \frac{1}{t}\int_0^{t} f(r^{i,{N_t}}_{\kappa_n(s)}) ds\right) \right)^2 \right]^{1/2}.		
\end{align*}
The mean-square error decomposition reads 
\begin{align*}
	&\textit{mse}(\cAesaea(t,n,{N_t})) \\
& \qquad  \lesssim     \left| \int_{\bR^d} f(x) (\pi (dx) - \mu_{t}(dx))  \right|    +   \left| \bE[ f(y_{t}) ]  - \bE[f(z^{1,{N_t}}_{t})]  \right|   \\
&\qquad   +  \left| \bE[  f(z^{1,{N_t}}_{t})]    - \bE[ f(r^{1,{N_t}}_t)]  \right|  
 +  \left(\bE \left[\left(\bE[ f(r^{1,{N_t}}_t)]   - \frac{1}{{N_t}} \sum_{i=1}^{N_t} \left( \frac{1}{t}\int_0^{t} f(r^{i,{N_t}}_{\kappa_n(s)}) ds\right) \right)^2 \right] \right)^{1/2 }
\end{align*}
Reasoning as before, with \eqref{as weak particle error} replaced by \eqref{as weak ergodic particle error}, we have
\[
\textit{mse}(\cAesaea(t,n,{N_t}) )\lesssim e^{-\lambda t} + e^{-\lambda t}({N_t})^{-1} + 1/n+ 1/\sqrt{t{N_t}}.
\]
Note that because of the variance (last term) there is no benefit of exponential decay in  $t$ assumed in \ref{as weak ergodic particle error}. Again we fix $\epsilon$. The following choice of the parameters ensures that $\textit{mse}(\cAesaea(t,n,{N_t})) \lesssim~\epsilon$:
\[
t \approx \lambda^{-1} \log(\epsilon^{-1}), \quad t{N_t} \approx \epsilon^{-2}, \quad n \approx \epsilon^{-1}. 
\]
Notice that the choice implies that $t N_t$ to be "constant". Hence we chose $N_t= N t^{-1}$. Which in the case of $t \approx \lambda^{-1} \log(\epsilon^{-1})$ implies that $N\approx \epsilon^{-2}/\log(\epsilon^{-1}) $.  

Now we study the computational cost of simulating self-interacting diffusions (for the nonlinear interacting kernel). Note that, due nonlinear interactions at every step of the Euler scheme, we have $N_t$ particles and each particle interacts with itself from all the past times-steps. Recall also that we take $n$ time steps in each unit time interval so that overall number of steps on the interval $[0,t]$ is $tn$. With that in mind and the fact that we take $N_t= N t^{-1}$ we have
\begin{align*}
&\textit{cost}({\cAesaea(t,n,N) })  
= N_{1/n}^2 1 + N_{2/n}^2 2 + N_{3/n}^2 3 + \ldots	+ N_{tn/n}^2 tn \\
&=  (nN)^2 (1+1/2 + 1/3 + \ldots + 1/tn) = (nN)^2 \sum_{k=1}^{tn} 1/k  \approx (nN)^2 \log(1+ tn). 
\end{align*}
Hence the cost for the set of parameters $\{t \approx \lambda^{-1} \log(\epsilon^{-1}),~N\approx \epsilon^{-2}/\log(\epsilon^{-1}),~ n \approx \epsilon^{-1}\}$ is
\[
\textit{cost}({\cAesaea(t,n,N) })  \approx \epsilon^{-6}
\] 

Note also that one can also take 
\[
t \approx \epsilon^{-2}, \quad (N_t =1 \implies t N_t = \epsilon^{-2}),  \quad n \approx \epsilon^{-1}
\]
to ensure $\textit{mse}(\cAesaea(t,n,{N_t})) \lesssim \epsilon$.
This second choice, with $N_t$ not varying with time, leads to 
\begin{align*}
&\textit{cost}({\cAesaea(t,n,N) })  \\
&= N^2 1 + N^2 2 + N^2 3 + \ldots	+ N^2 tn =  N^2 \frac{1}{2} tn (1 + tn)
%\approx (N tn)^2
 \approx \epsilon^{-6}.
 \end{align*}

%\ls{Like I said maybe there are better ways of decomposing the error} 

Let us consider ensemble implementation of the above algorithm. Reasoning as before, by the Assumption \ref{as weak ergodic particle error}, we have
\[
\textit{mse}(\cAesaea(t,n,{N_t}) )\lesssim e^{-\lambda t} + e^{-\lambda t}({N_t})^{-1} + 1/n+ 1/\sqrt{t{N_t}M}.
\]
To balance the first two terms on the right hand side we take $N_t=1$. With that choice we then chose $M$ so that $e^{-\lambda t}= (t M)^{-1/2}$, i.e 
$M=e^{2 \lambda t } t^{-1}$. With this choices we have
\[
\textit{mse}(\cAesaea(t,n,{N_t}) )\lesssim e^{-\lambda t} + 1/n,
\]
to make that error to be less then $\epsilon$ we take $t=\lambda^{-1} \log{\epsilon^{-1}}$ and  $n=\epsilon^{-1}$. Note that this leads to
$M=e^{\log(\epsilon^{-2})}(\log(\epsilon^{-1}))^{-1}=\epsilon^{-2}(\log(\epsilon^{-1}))^{-1}$. Reasoning as before the cost with $N_t=1$ is
\begin{align*}
&\textit{cost}({\cAcsaea(t,n,N_t,M) })  
\approx M (tn)^2 \approx \epsilon^{-4}\log(\epsilon^{-1}). 
\end{align*}
Hence one more time we achieved order better computational cost in comparison to 
naive estimator. Note that the presented analysis implies that ensemble of $M$ independent self-interacting diffusions yields best result.

\section{Conclusion and perspectives }\label{conc}
We presented a number of different algorithms for the approximation of the invariant measure of the McKean-Vlasov SDE. We have achieved one order better computational costs comparing to naive particle based estimator. On algorithmic side possible extensions are that may consider fixed length window for self-interacting diffusion and it also possible to study Multilevel Monte Carlo strategies in this setup. Overall we anticipate that it will be possible to bring the cost of simulating particle system to the same level as for standard independent copies of SDEs.

\bibliographystyle{alpha}%wmaainf}%alpha}%plain}
\bibliography{Particles}  % ``name``.bib is the name of thedatabase
% \bib, bibdiv, biblist are defined by the amsrefs package.
%%%%%%%%%%%%%%%%%%%%%%%%%%%%%%%%%%%%%%%%%%%%%%%%%%%%%%%%%%%%%%%%%%%%%%%%%%
%%%% \END BIBLIOGRAPHY
%%%%%%%%%%%%%%%%%%%%%%%%%%%%%%%%%%%%%%%%%%%%%%%%%%%%%%%%%%%%%%%%%%%%%%%%%%
\end{document}